\documentclass[11pt,reqno,tbtags]{amsart}
\usepackage{amssymb}
\numberwithin{equation}{section}

\newtheorem{theorem}{Theorem}[section]
\newtheorem{lemma}[theorem]{Lemma}

\newtheorem{corollary}[theorem]{Corollary}

\theoremstyle{definition}

\newtheorem{remark}[theorem]{Remark}

\newtheorem*{acks}{Acknowledgements}

\newenvironment{thmxenumerate}{\begin{enumerate}
\setlength{\leftmargin}{0pt}
\setlength{\itemindent}{0pt}
 }{\end{enumerate}}

\newcounter{thmenumerate}
\newenvironment{thmenumerate}
{\setcounter{thmenumerate}{0}%
 \def\item{\ifnum\thethmenumerate=0\else\par\fi 
 \addtocounter{thmenumerate}{1}\textup{(\roman{thmenumerate})\enspace}}}
{}
\newcommand\pfitem[1]{\par(#1):}

\newcounter{xenumerate}
\newenvironment{xenumerate}
  {\begin{list}
    {\upshape(\roman{xenumerate})}
    {\setlength{\leftmargin}{0pt}
     \setlength{\rightmargin}{0pt}
     \setlength{\labelwidth}{0pt}
     \setlength{\itemindent}{\labelsep}
     \setlength{\topsep}{0pt}
     \usecounter{xenumerate}} }
  {\end{list}}

\newcommand{\refT}[1]{Theorem~\ref{#1}}
\newcommand{\refC}[1]{Corollary~\ref{#1}}
\newcommand{\refL}[1]{Lemma~\ref{#1}}

\newcommand{\refS}[1]{Section~\ref{#1}}
\newcommand{\refTab}[1]{Table~\ref{#1}}
\newcommand{\refand}[2]{\ref{#1} and~\ref{#2}}

\begingroup
  \divide\count255 by 60
  \multiply\count255 by -60
  \advance\count255 by \time
  \ifnum \count255 < 10 \xdef\klockan{\the\count1.0\the\count255}
  \else\xdef\klockan{\the\count1.\the\count255}\fi
\endgroup

\newcommand\set[1]{\ensuremath{\{#1\}}}

\newcommand\xpar[1]{(#1)}
\newcommand\bigpar[1]{\bigl(#1\bigr)}
\newcommand\Bigpar[1]{\Bigl(#1\Bigr)}

\newcommand\bigsqpar[1]{\bigl[#1\bigr]}
\newcommand\Bigsqpar[1]{\Bigl[#1\Bigr]}

\def\rompar(#1){\textup(#1\textup)}    

\newcommand\parfrac[2]{\Bigpar{\frac{#1}{#2}}}

\def\xexp(#1){e^{#1}}

\newcommand\mtoo{\ensuremath{{m\to\infty}}}
\newcommand\ntoo{{\ensuremath{n\to\infty}}}

\newcommand\bmin{\wedge}

\newcommand\iid{i.i.d.\spacefactor=1000}    
\newcommand\ie{i.e.\spacefactor=1000}
\newcommand\eg{e.g.\spacefactor=1000}

\newcommand\cf{cf.\spacefactor=1000}

\newcommand\whp{\textbf{whp}}

\newcommand{\tend}{\longrightarrow}
\newcommand\dto{\overset{\mathrm{d}}{\tend}}
\newcommand\pto{\overset{\mathrm{p}}{\tend}}

\newcommand\ucpto{\overset{\mathrm{ucp}}{\tend}}

\newcommand\bbN{\ensuremath{\mathbb N}}

\newcommand\E{\operatorname{\mathbb E{}}}
\renewcommand\P{\operatorname{\mathbb P{}}}
\newcommand\Var{\operatorname{Var}}

\newcommand\Exp{\operatorname{Exp}}
\newcommand\Po{\operatorname{Po}}

\newcommand\Ge{\operatorname{Ge}}

\def\[#1]{[\![#1]\!]}

\newcommand\ind[1]{[\![#1]\!]}

\newcommand\ga{\alpha}

\newcommand\gd{\delta}

\newcommand\gf{\varphi}

\newcommand\eps{\varepsilon}

\renewcommand\phi{\gf $\marginpar{SMISK!}$ }

\newcommand\cT{{\mathcal T}}

\newcommand\qi{^{-1}}

\renewcommand{\=}{:=}

\newcommand\intot{\int_0^t}
\newcommand\intoa{\int_{0}^{\ga}}

\newcommand\et{e^{-t}}
\newcommand\eet{1-e^{-t}}

\newcommand\ctmn{\ensuremath{\cT_{m,n}}}

\newcommand\ctt{\ensuremath{\cT(t)}}

\newcommand\dxit{\ensuremath{d_i^\Xi(\cT)}}
\newcommand\dxIt{\ensuremath{d_I^\Xi(\cT)}}

\newcommand\nxkt{\ensuremath{n_k^\Xi(\cT)}}

\newcommand\dxxt{\ensuremath{d^\Xi(\cT)}}
\newcommand\da[1]{\ensuremath{D_\ga^{#1}}}
\newcommand\dax{\ensuremath{\da\Xi}}
\newcommand\dau{\ensuremath{\da\U}}
\newcommand\dal{\ensuremath{\da\li}}
\newcommand\dae{\ensuremath{\da\ei}}

\newcommand\pax{\ensuremath{p^\Xi_\ga}}
\newcommand\pau{\ensuremath{p^\U_\ga}}
\newcommand\pal{\ensuremath{p^\li_\ga}}
\newcommand\pae{\ensuremath{p^\ei_\ga}}

\newcommand\dtmn[1]{\ensuremath{d^{#1}(\ctmn)}}
\newcommand\dxtmn{\dtmn{\Xi}} 

\newcommand\dutmn{\ensuremath{d^\U(\ctmn)}}
\newcommand\dltmn{\ensuremath{d^\li(\ctmn)}}
\newcommand\detmn{\ensuremath{d^\ei(\ctmn)}}

\newcommand{\duj}{d_j^\U}
\newcommand{\tdu}{\widetilde d^\U}
\newcommand{\tdau}{\widetilde D_\ga^\U}
\newcommand{\nxk}{n_k^\Xi}

\newcommand{\nuk}{n_k^\U}

\newcommand{\nek}{n_k^\ei}
\newcommand{\nlk}{n_k^\li}

\newcommand\duxt[1]{\ensuremath{d_{#1}^\U(\cT)}}
\newcommand\dujt{\duxt{j}}

\newcommand\dut{\duxt{}}

\newcommand\gaix{\frac1{\ga}}
\newcommand\gai{\ga\qi}

\newcommand{\ei}{\textup{\textsf{E}}}
\newcommand{\li}{\textup{\textsf{L}}}
\newcommand{\U}{\textup{\textsf{U}}}

\newcommand{\el}{\set{\li,\ei}}
\newcommand{\elu}{\set{\li,\ei,\U}}
\newcommand\dd{\,d}

\newcommand{\pgf}{probability generating function}

\newcommand\lhs{left hand side}

\newcommand\nc{\newcommand}
\newcommand\eu{e^{-u}}
\nc{\clt}{C_\ell(t)}
\nc{\cklt}{\sum_{\ell\ge k} \clt}
\nc\fk[1]{f^{(k)}_{#1}}

\newcommand{\maple}{\texttt{Maple}}

\newcommand\REM[1]{\texttt{[#1]}\marginal{XXX}}

\begin{document}
\title
{Individual displacements in hashing with coalesced chains}

\date{February 10, 2005} 

\author{Svante Janson}
\address{Department of Mathematics, Uppsala University, PO Box 480,
S-751 06 Uppsala, Sweden}
\email{svante.janson@math.uu.se}
\urladdr{http://www.math.uu.se/\~{}svante/}


\begin{abstract} 
We study the asymptotic distribution of the displacements in
hashing with coalesced chains, for both late-insertion and
early-inser\-tion. Asymptotic formulas for means and variances follow.
The method uses Poissonization and some stochastic calculus.
\end{abstract}

\maketitle


\section{Introduction}\label{S:Intro}
The standard version of hashing with coalesced chains,
due to Williams \cite{Williams}
can be
described as follows, where $n$
and $m$  are integers with $0\le n\le m$. 
(See further Knuth \cite[Section 6.4, in particular Algorithm 6.4.C]{KnuthIII}
and the monograph by Vitter and Chen \cite{VitterC}.)

\begin{quote}
$n$ items $x_1,\dots,x_n$ are placed sequentially into a table with
$m$ cells $1,\dots,m$, using $n$ integers $h_i\in\set{1,\dots,m}$.
Each cell contains a link field, initially null.
Item $x_i$ is inserted into cell $h_i$ if it is empty; otherwise we
follow the links from cell $h_i$ until we reach the end of the chain
(signalled by a null link), we add a link to an empty cell (which is
chosen as the empty cell with largest index) and store the item there.
\end{quote}

For our probabilistic treatment, we assume that 
each of the $m^n$ possible hash sequences $(h_i)_1^n$ is equally likely;
in other words, 
the hash addresses $h_i$ are independent random numbers, uniformly
distributed on $\set{1,\dots,m}$.

The \emph{displacement} $d_i$ of an item $x_i$ is the number of links we
have to follow from $h_i$ until we find $x_i$. Large displacements
make both insertion and searching less efficient, so it is desirable
to keep the displacements small.
(Two different but related quantities are used in other papers
to measure the efficiency:
\emph{The number of probes} to find the item $x_i$ in the table is 
$d_i+1$.
\emph{The number of key comparisons} to find the item is also
$d_i+1$. 
This should be noted when comparing the results below with
other papers.)

The items are thus arranged in linked chains in the hash table. If a
new item hashes to an empty cell, a new chain with that single item is
created. If a new item hashes to a cell in an existing chain, then
that chain grows by addition of a formerly empty cell.
It was shown by Chen and Vitter \cite{ChenVitter83} 
and Knott \cite{Knott}
that the average
performance could be improved by modifying the algorithm above,
inserting the new item at a different place in an existing chain.
We will therefore study two versions of hashing with coalesced chains:
\begin{itemize}
\item[\li] \emph{Late-insertion} (LISCH).
The standard version described above where the new item is inserted last
in its chain.

\item[\ei]
\emph{Early-insertion} (EISCH)
\cite{ChenVitter83}, \cite{Knott}, \cite{VitterC}.
If cell $h_i$ is occupied, item $x_i$ is inserted into an empty cell
as above, but this cell is linked into the chain immediately after
$h_i$.
(I.e., if the first free cell is $j$ and the link from $h_i$ points to
$k$ (null or not), then this link is reset to $j$, and the link field
in $j$ is set to $k$.)
This method gives the smallest average displacement among all possible
insertion schemes 
\cite[Theorem 5.2]{VitterC}.
\end{itemize}

Note that the insertion of a sequence of items
results in the same set of occupied cells in both versions,
and that this set is partitioned into chains in the same way,
but that the order in the chains, and thus the individual
displacements, 
may differ.

Our main result is \refT{T1} below (together with its refinement
\refT{T2}), which gives the asymptotic distribution of the
displacements in a random hash table under both insertion methods;
we consider also the case of unsuccessful searches.
As corollaries we easily find earlier known asymptotic formulas for
the average displacements and for the variances of them (some of the latter
may be new). These asymptotic distributions are studied further in
\refS{Sdistr}; some numerical values are given in \refTab{Ta1}.
The proofs are given in \refS{Sproofs}. They are based on
Poissonization, regarding the items as arriving at random times.

\begin{remark}
An interesting variation of the algorithm above 
\cite[Exercise 6.4-43]{KnuthIII}, 
discussed in detail
by Vitter and Chen
\cite{VitterC}, 
is to choose a number $m_1<m$  and reserve the
last $m-m_1$ cells as a ``cellar''
for the undisturbed growth of the chains. We then assume
that the hash addresses $h_i\in\set{1,\dots,m_1}$, and use exactly
the same algorithms as above. (These versions are called LICH and EICH
in \cite{VitterC}.)
It is shown in \cite{VitterC}
that for given $n$ and $m$, a suitable choice of $m_1$ will improve
the average performance. 

In this setting, it is also interesting to consider a third version
\emph{varied-insertion} 
(VICH) \cite{VitterC}, which behaves like $\ei$
except that when the chain from the hash address contains a cellar
cell, the new item is inserted after the last cellar cell.
It is shown in \cite[Chapter 5]{VitterC}
that this method gives the minimum average among all insertion methods
satisfying a weak assumption.

We have not yet investigated the versions with cellar in detail, but
it seems that 
our methods could be used with some additional work to find the
asymptotic distributions of the displacements in these cases too. (The
means are given in \cite{VitterC}.) 
\end{remark}

\begin{remark}
Corresponding results for hashing with linear probing are given by
Janson \cite{SJ157} and Viola \cite{Viola03}. Note that, as remarked in
\cite{KnuthIII},  the average displacement for 
linear probing tends to infinity if $n,m\to\infty$ with $n/m\to1$,
while for the chained hashing studied here, it stays bounded also in
the extreme case $n=m$ of a full table.
\end{remark}

\begin{acks}
This research was begun on the transatlantic flight to the Knuthfest
at Stanford in honour of 
Donald Knuth's 64th birthday in January 2002, where some of the
results were presented, and continued on another flight to the 
Analysis of Algorithms workshop at MSRI, Berkeley, June 2004;
I thank KLM for providing long and quiet periods of work.
I further thank Donald Knuth for helpful remarks.
\end{acks}

\section{Notation and results}\label{S:prel}
By a hash table $\cT$ we mean not only the final table, but also its
construction history; moreover, we consider the two version above together.
Formally, a hash table can be regarded as encoded by the numbers $m$
and $n$ and the sequence
$(h_1,\dots,h_n)$ of hash addresses.

Our prime object of study is the random hash table \ctmn{} with $m$
cells and $n$ items ($0\le n\le m$) and the hash addresses
$h_1,\dots,h_n$ \iid{} random variables, uniformly distributed on
\set{1,\dots,m}.

We denote the two insertion policies defined in the introduction by
\li{} and \ei, and use $\Xi$ to denote any of these.

Given a hash table $\cT$ (with $m$ cells and $n$ items), random or
not, and a policy $\Xi\in\el$, we 
let \dxit{} be the (final) displacement of the $i$:th item, 
$1\le i\le n$, and
\begin{equation*}
\nxkt\=\#\set{i:\dxit=k},\qquad k=0,1,\dots,
\end{equation*}
the number of items with displacement $k$.
Note that 
\begin{equation}
\label{sumnx}
\sum_k\nxkt=n. 
\end{equation}

If $n>0$, we let \dxxt{} denote a randomly chosen displacement in a given hash
table $\cT$ using policy $\Xi$, \ie{} the random variable 
\dxIt{} where $I\in\set{1,\dots,n}$ is a random index with a uniform
distribution.
Thus, given $\cT$, \dxxt{} has the distribution
\begin{equation}\label{pt}
\P\bigpar{\dxxt=k\mid\cT}=\frac1n \nxkt.
\end{equation}

Similarly,
we let \dujt{} denote the number of occupied
cells encountered in an \emph{unsuccessful} search starting at hash
address $j$, $1\le j\le m$,
and let \dut{} denote the number of occupied cells
encountered in a random unsuccessful  search,
\ie{}
$\dut\=\duxt{J}$, where $J\in\set{1,\dots,m}$ is a
uniformly distributed random index.
(Note that in an unsuccessful search starting at $j$,
the number of key comparisons
equals $d_j^\U$, while the number of
probes, $\tdu_j$ say, 
is $d_j^\U$ if $d_j^\U\ge1$ and 1 if $d_j^\U=0$; 
\ie{}, $\tdu_j\=\max(d_j^\U,1)$).)
We further let
\begin{equation*}
\nuk(\cT)\=\#\set{j:\dujt=k},\qquad k=0,1,\dots,
\end{equation*}
and note that now, in contrast to \eqref{sumnx},
\begin{equation}
\label{sumnu}
\sum_k\nuk(\cT)=m.
\end{equation}
Thus, \cf{} \eqref{pt}, given $\cT$, \dut{} has the distribution
\begin{equation*}\label{ut}
\P\bigpar{\dut=k\mid\cT}=\frac1m \nuk(\cT).
\end{equation*}

Our main results is the following theorem, giving the
asymptotic distribution of the 
displacement of a random item in a random hash table $\ctmn$,
together with the corresponding quantity for an unsuccessful search.

\begin{theorem}\label{T1}
Suppose that $m,n\to\infty$ with $0< n\le m$ and $n/m\to\ga\in[0,1]$.
Then, for every $k=0,1,\dots$,
\begin{xenumerate}
\item 
\begin{multline*}
\P\bigpar{\dtmn\U=k}
=\frac1m \E\bigpar{\nuk(\ctmn)}
\\
\to \pau(k)
\=
\begin{cases}
  1-\ga, & k=0, 
\\
\intoa(1-\ga+t)\bigpar{\eet}^{k-1}\dd t, & k\ge1;
\end{cases}	
\end{multline*}
\item 
\begin{multline*}
\P\bigpar{\dtmn\li=k}
=\frac1n \E\bigpar{\nlk(\ctmn)}
\\
\to \pal(k) 
\= 
\begin{cases}
  1-\ga/2, & k=0, 
\\
\gaix\intoa\bigpar{\ga-t-(\ga-t)^2/2}\bigpar{\eet}^{k-1}\dd t, & k\ge1;
\end{cases}	
\end{multline*}
\item 
\begin{multline*}
\P\bigpar{\dtmn\ei=k}
=\frac1n \E\bigpar{\nek(\ctmn)}
\\
\to \pae(k)
\=
\begin{cases}
  1-\ga/2, & k=0, 
\\
\gaix\intoa (\ga-t)\et\bigpar{\eet}^{k-1}\dd t, & k\ge1.
\end{cases}	
\end{multline*}
\end{xenumerate}
(For $\ga=0$, $p_0^\li(k)=p_0^\ei(k)=0$ when $k\ge1$.)
For every $\Xi\in\elu$ and $\ga\in[0,1]$,
$\set{\pax(k)}_{k=0}^\infty$ is a probability distribution on
$\bbN$. 
If $\dax$ is a random variable with this distribution,
\ie{}
$\P(\dax=k)=\pax(k)$, then these results can be written
\begin{equation}\label{t1}
\dxtmn\dto \dax.
\end{equation}
Moreover, all  moments converge in \eqref{t1}, \ie,
$\E\bigpar{\dxtmn}^r\to \E\bigpar{\dax}^r$
for every $r\ge0$.
\end{theorem}

It follows immediately that for the number of probes in an
unsuccessful search, we have 
\begin{equation}
  \label{jesper}
\tdu(\ctmn)\dto\tdau\=\max(\dau,1),
\end{equation}
again with convergence of all moments.

As a corollary, we find the asymptotics for the expectations; these
have earlier been derived, 
together with exact formulas for
$\E\dxtmn$,
by Knuth \cite{KnuthIII} and Vitter and Chen
\cite{VitterC}
(in equivalent forms for the number of probes or key comparisons).

\begin{corollary}\label{C1}
Suppose that $m,n\to\infty$ with $0< n\le m$ and $n/m\to\ga\in[0,1]$.
Then
\begin{align*}
\E(\dutmn)& \to \E \dau=\frac14\bigpar{e^{2\ga}-1}+\frac{\ga}2, 
\\
\E(\dltmn)& \to \E \dal=\frac1{8\ga}\bigpar{e^{2\ga}-1}
+\frac{\ga}4-\frac14, 
\\
\E(\detmn)& \to \E \dae=\frac1\ga\bigpar{e^{\ga}-1-\ga}. 
\end{align*}
\end{corollary}

\begin{remark}
Note that the expected number of probes is $\E(\dltmn)+1$ or
$\E(\detmn)+1$
for a successful search, and, \cf{} \eqref{jesper}, 
\begin{align*}
\E(\tdu(\ctmn))=\E(\dutmn)+\frac{m-n}m 
\to \E(\tdau)=\E(\dau)+1-\ga
\end{align*}
for an  unsuccessful search. 
\end{remark}

\refT{T1} similarly yields asymptotic formulas for higher moments too;
in particular we have the following results for the variance.

\begin{corollary}\label{C2}
Suppose that $m,n\to\infty$ with $0< n\le m$ and $n/m\to\ga\in[0,1]$.
Then
\begin{align*}
\Var(\dutmn)& \to
\Var\dau
\\&\hskip2em
=
  -\tfrac1{16}e^{4\ga}+\tfrac49e^{3\ga}-\bigpar{\tfrac14\ga+\tfrac18}e^{2\ga} 
-\tfrac14\ga^2+\tfrac{5}{12}\ga-\tfrac{37}{144},
\\[3pt]
\Var(\dltmn)& \to
\Var\dal
=
-\tfrac1{64}\parfrac{e^{2a}-1}{a}^2
+\frac{64e^{2a}+37e^a+37}{432}\cdot\frac{e^a-1}{a}
\\&\hskip12em
-\tfrac1{16}e^{2a}
-\tfrac1{16}a^2
+\tfrac{5}{24}a
-\tfrac7{36},
\\[3pt]
\Var(\detmn)& \to
\Var\dae
=
\frac{\ga-2}{2}\parfrac{e^\ga-1}{\ga}^2+2\frac{e^\ga-1}{\ga}-1,
\\
\intertext{and, for the number of probes in an unsuccessful search,}
\Var(\tdu(\ctmn))
&\to\Var(\tdau)
\\&\hskip2em
=
-\tfrac1{16}e^{4\ga}+\tfrac49e^{3\ga}+\bigpar{\tfrac14\ga-\tfrac58}e^{2\ga} 
-\tfrac14\ga^2-\tfrac{1}{12}\ga+\tfrac{35}{144}.
\end{align*}
\end{corollary}

The asymptotic formula for $\Var(\tdu(\ctmn))$, together with an exact
formula, is given in Knuth \cite[Answer 6.4-40]{KnuthIII}
and in Vitter and Chen \cite{VitterC}; 
the corresponding results for $\Var(\dutmn)$ follow easily.
(The numerical result in \cite[Answer 6.4-40]{KnuthIII} and
\cite{VitterC} for the case $\ga=1$, when $\tdau=\dau$, should be $2.65$.)
The asymptotics of $\Var(\dltmn)$ are given in \cite{VitterC}.
We do not know whether the asymptotic of $\Var(\detmn)$ have been
published earlier.

Consider a computer program where a large
hash table is constructed once, and then used many times
for finding the items. 
We assume that each item in the table is equally
likely to be requested, and that each choice is independent of the
previous ones.
We therefore have two levels of randomness: First we construct a
random hash table $\cT$ with some 
displacements $(d_i)$. Keeping $\cT$ fixed and 
choosing a random index $I\in \set{1,\dots,n}$,
we obtain the random displacement
$d(\cT)=d_I$.
As the program runs with many searches in the hash table, 
the search times then are (functions of) independent observations of
this random 
variable. It is thus interesting to study the distribution of this
random variable and its properties such as its mean and variance.
Note that this distribution depends on the hash table $\cT$,
which is itself random;
another run of the program yields another $\cT$ and another set of
displacements.
Hence the distribution of the displacement $d(\cT)$ 
is a random distribution and its mean
$\E(d(\cT)|\cT)$ 
and variance $\Var(d(\cT)|\cT)=\E(d(\cT)^2|\cT)-\E(d(\cT)|\cT)^2$ 
are random variables.
In other words, we study the conditional distribution of $d(\cT)$
given $\cT$.

We can refine the results above by conditioning on $\ctmn$.
The following theorem says that we still have the same limits, now
with convergence in probability.
In other words, 
different realizations of $\ctmn$ have (with large probability) almost the same
distribution of the displacements, so
a typical instance of the random hash table
$\ctmn$ has its displacements distributed as the average studied in \refT{T1}.

\begin{theorem}
  \label{T2}
Suppose that $m,n\to\infty$ with $0< n\le m$ and $n/m\to\ga\in[0,1]$.
Then, for every $k=0,1,\dots$, with $\pax(k)$ defined in \refT{T1},
\begin{thmxenumerate}
  \item
$\displaystyle
\P\bigpar{\dutmn=k\bigm|\ctmn}
=m\qi \nuk(\ctmn)
\pto \pau(k)$,
\item
$\displaystyle
\P\bigpar{\dltmn=k\bigm|\ctmn}
=n\qi \nlk(\ctmn)
\pto \pal(k)$,
\item
$\displaystyle
\P\bigpar{\detmn=k\bigm|\ctmn}
=n\qi \nek(\ctmn)
\pto \pae(k)$.
\end{thmxenumerate}
\end{theorem}

\begin{remark}
A more fancy formulation of \refT{T2} is that the distribution of
$\dxtmn$ converges to $\pax$ in probability, in the space of all
probability measures on $\bbN$, equipped with the weak topology
(which coincides with the $\ell^1$ topology on this space); see
\cite{Bill} for definitions. 
\end{remark}

Moment convergence holds in \refT{T2} too, \ie{} conditioned on
$\ctmn$. 

\begin{theorem}
  \label{T3}
Suppose that $m,n\to\infty$ with $0< n\le m$ and $n/m\to\ga\in[0,1]$.
Then, for every $r\ge0$ and $\Xi\in\elu$,
\begin{equation*}
  \E\bigpar{\dxtmn^r\mid\ctmn}\pto \E\bigpar{\dax}^r.
\end{equation*}
In particular, the conditional mean and variance, given the
hash table, converge in probability to the limits in Corollaries
\refand{C1}{C2}. 
\end{theorem}

\section{The asymptotic distributions}\label{Sdistr}

We give some further results on the probability distributions
$\pax(k)$ defined in \refT{T1}; we assume $\ga>0$.
We omit the proofs.
(Several of the results below were obtained with the help of \maple.)

It follows directly from the definitions in \refT{T1} that
\begin{equation*}
\pax(k)=O\bigpar{1-e^{-\ga}}^k=O\bigpar{1-e^{-1}}^k;   
\end{equation*}
hence the 
probabilities decrease geometrically. More refined asymptotics can
easily be derived (we omit the details); we have, as $k\to\infty$,
for $\ga>0$,
\begin{align*}
\pau(k) &\sim e^\ga k^{-1} \bigpar{1-e^{-\ga}}^k,
\\
\pau(k) &\sim \gai e^\ga(e^\ga-1) k^{-2} \bigpar{1-e^{-\ga}}^k,
\\
\pau(k) &\sim \gai(e^\ga-1)k^{-2} \bigpar{1-e^{-\ga}}^k.
\end{align*}
In particular, the probability of an extremely large displacement is
about $e^\ga$ as large for late-insertion as for early-insertion.

The \pgf{s} for $\dax$ follow also easily from the formulas in \refT{T1}:
\begin{align*}
\E x^{\dau}
&=
\sum_{k=0}^\infty \pau(k)x^k
=
1-\ga+
x\intoa\frac{1-\ga+t}{1-x+x\et}\dd t,
\\
\E x^{\dal}
&=
\sum_{k=0}^\infty \pal(k)x^k
=
1-\frac\ga2+
\frac x\ga\intoa\frac{\ga-t-(\ga-t)^2/2}{1-x+x\et}\dd t,
\\
\E x^{\dae}
&=
\sum_{k=0}^\infty \pae(k)x^k
=
1-\frac\ga2+
\frac x\ga\intoa \frac{\ga-t} {(1-x)e^t+x}\dd t.
\end{align*}
These integrals can be evaluated in terms of the dilog
function (and for $\li$ also polylog), 
but we do not know any simple
form.
The generating functions are analytic for 
$|x|<r(\ga)\=\bigpar{1-e^{-\ga}}\qi$,
with a singularity at $r(\ga)$. 

The integrals defining $\pax(k)$ are easily evaluated for small $k$.
We find, for example, 
\begin{align*}
\pau(1)&=\ga-\tfrac12\ga^2, 
&
\pau(2)&=2e^{-a}-2+2\ga-\tfrac12\ga^2,
\\
\pal(1)&=\tfrac12\ga-\tfrac16\ga^2,
&
\pal(2)&=2\frac{1-e^{-\ga}}{\ga}-2+\ga-\tfrac16\ga^2,
\\
\pae(1)&=\frac{e^{-\ga}-1+\ga}\ga,
&
\pae(2)&=\frac12-\frac{1-e^{-\ga}}{\ga}+\frac{1-e^{-2\ga}}{4\ga}.
\end{align*}
No simple pattern is seen, and we leave further investigation to the
reader.

Numerical values for $\pax(0),\dots,\pax(10)$, the tail
$\sum_{11}^\infty\pax(k)$, the mean $\E\dax$ and the variance
$\Var\dax$ are given for $\ga=0.5$ and 1 (half-full and full tables)
in \refTab{Ta1}.

\begin{table}[htb]
  \begin{tabular}{ r| l | l | l | l | l | l }
$k$ & $p_{0.5}^\U(k)$ & $p_{0.5}^\li(k)$ & $p_{0.5}^\ei(k)$
    & $p_{1}^\U(k)$ & $p_{1}^\li(k)$ & $p_{1}^\ei(k)$ \\
\hline 
0& 0.5& 0.75& 0.75& 0.0& 0.5& 0.5\\
1& 0.375& 0.2083& 0.2130& 0.5& 0.3333& 0.3679\\
2& 0.0881& 0.0322& 0.0291& 0.2358& 0.0976& 0.0840\\
3& 0.0252& 0.0070& 0.0059& 0.1200& 0.0376& 0.0280\\
4& 0.0078& 0.0018& 0.0014& 0.0638& 0.0163& 0.0110\\
5& 0.0026& 0.00049& 0.00038& 0.0349& 0.0076& 0.0048\\
6& 0.00086& 0.00014& 0.00011& 0.0194& 0.0037& 0.0022\\
7& 0.00030& 0.000043& 0.000032& 0.0110& 0.0019& 0.0011\\
8& 0.00010& 0.000014& 0.000010& 0.0063& 0.0010& 0.0005\\
9& 0.00004& 0.000004& 0.000003& 0.0036& 0.0005& 0.0003\\
10& 0.00001& 0.000001& 0.000001& 0.0021& 0.0003& 0.0001\\
\hline
$\ge11$ &
0.000007& 0.0000007& 0.0000005& 0.0031& 0.0003& 0.0002\\ 
\hline
\hline
$\E$& 0.6796& 0.3046& 0.2974& 2.0973& 0.7986& 0.7183\\
$\Var$& 0.7394& 0.3565& 0.3324& 2.6533& 1.2799& 0.9603\\
\hline
  \end{tabular}
\caption{Some numerical values\strut}
\label{Ta1}
\end{table}

\section{Proofs}\label{Sproofs}

To prove the theorems, 
we randomize the times the items are inserted in the table by
Poissonization: We assume that items with hash address $i$ arrive
according to a Poisson process with intensity $1$, the $m$ different
Poisson processes being independent. We let $\ctt$ denote the hash
table at time $t$, when there are $\Po(t)$ items with each hash address.
For simplicity, we write 
$\nxk(t)\=\nxk(\ctt)$.

Combining the $m$ individual Poisson processes, we see that 
the items $x_1,x_2,\dots$
arrive according to a Poisson process with intensity
$m$; we call the arrival times $\tau_1,\tau_2,\dots$ (we may assume
that these are distinct).
We really have to stop at $\tau_m$, since the table then is full, but
it is convenient to think of the hashing as continuing for ever, with
the chains growing into a virtual, infinitely large attic; no new chains are
created after $\tau_m$. 

The hash addresses of the items $x_1,x_2,\dots$
are independent and uniformly distributed,
so except for the random time scale, this is the situation we want to
study. More precisely, $\cT(\tau_n)=\ctmn$ for $0\le n\le m$.

Note that $\tau_m\approx1$; more precisely, $\tau_m/m\pto1$ as \mtoo,
as shown in \refL{Ltau} below.

We will consider stochastic processes defined on $[0,\infty)$
(although we mainly are interested in $0\le t\le 1$).
We say that such a process $X(t)$ is increasing if $X(s)\le X(t)$
whenever $s\le t$. We let $\ucpto$ denote convergence 
\emph{uniformly on compacts in probability} (ucp), \ie{} $X_n\ucpto X$ if
$\sup_{0\le t\le u}|X_n(t)-X(t)|\pto0$ for every $u>0$. 

\begin{lemma}
  \label{LU}
Let, for each $n$, 
$X_n(t)$, $t\ge0$, be an increasing, stochastic
process, and let $f(t)$ be a continuous function on $[0,\infty)$.
If $X_n(t)\pto f(t)$ for every $t\ge0$, then $X_n(t)\ucpto f(t)$.
\end{lemma}
 
\begin{proof}
Fix $u>0$. Let $\eps>0$, and let $K$ be so large that if $\gd\=u/K$, then
$|f(s)-f(t)|<\eps$ if $s-t\le \gd$ and $0\le s\le t\le u$.
Since each $X_n$ is increasing, the limit $f(t)$ is too. Hence, 
if $(k-1)\gd\le t\le k\gd$,
\begin{multline*}
X_n((k-1)\gd)-f((k-1)\gd)-\eps 
\le X_n((k-1)\gd)-f(k\gd) 
\\
\le 
  X_n(t)-f(t) \le X_n(k\gd)-f((k-1)\gd) \le  X_n(k\gd)-f(k\gd)+\eps,
\end{multline*}
and, consequently,
\begin{equation*}
  \sup_{0\le t\le u}|X_n(t)-f(t)|
\le \sup_{0\le k\le K} |X_n(k\gd)-f(k\gd)|+\eps.
\end{equation*}
We know that $X_n(t)-f(t)\pto0$ for every $t\ge0$. We apply this for
$t=k\gd$, $k=0,\dots,K$, and find that \whp{} (\ie, with probability
$\to1$ as \ntoo) 
$|X_n(k\gd)-f(k\gd)|<\eps $ for $k=0,\dots,K$, and thus
$  \sup_{0\le t\le u}|X_n(t)-f(t)|<2\eps$.
Since $\eps>0$ is arbitrary, 
$\sup_{0\le t\le u}|X_n(t)-f(t)|\pto0$.
\end{proof}

Let $N(t)$ be the number of items that have arrived at time $t$.
Since each item is put into some empty cell, the number of empty cells
at time $t$ is $(m-N(t))_+$, \ie{} $m-N(t)$ for $t\le\tau_m$ and then 0.

\begin{lemma}
  \label{LN}
As \mtoo,
$ N(t)/m\ucpto t$ for $t\ge0$.
\end{lemma}
 
\begin{proof}
We have $N(t)\sim\Po(mt)$, and thus $N(t)/m\pto t$ as \mtoo{} for
every $t\ge0$. The convergence ucp follows from \refL{LU}.
\end{proof}

\begin{lemma}
  \label{Ltau}
Suppose that \mtoo{} and $n/m\to\ga$, with $0\le n\le m$.
Then $\tau_n\pto\ga$.
Consequently, if
$X_m\ucpto X$ for some stochastic processes $X_m$ and $X$, where $X(t)$
is continuous, then $X_m(\tau_n)\pto X(\ga)$.
\end{lemma}
\begin{proof}
$\tau_n$ is the sum of $n$ \iid{} waiting times, each $\Exp(1/m)$,
so $\E\tau_n=n/m\to\ga$ and $\Var\tau_n=n/m^2\to0$, whence
$\tau_n\pto\ga$ by Chebyshev's inequality.
(Alternatively,
this is a standard consequence of \refL{LN}:
For $\eps>0$, $\P(\tau_n>\ga+\eps) 
\le \P\bigpar{N(\ga+\eps)/m< n/m}
\to0$.
Similarly,
$\P(\tau_n<\ga-\eps) \to0$.)

The final assertion follows because \whp{} $\tau_n<\ga+1$, and then
$|X_m(\tau_n)-X(\ga)|
\le\sup_{s\le\ga+1} |X_m(s)-X(s)| + |X(\tau_n)-X(\ga)|\pto0$
\end{proof}

\subsection{Chains}
When an item arrives to an empty cell, a new chain of length 1 is
created. The chain then grows one unit each time it is hit. Hence each
chain, once created, grows according to a birth process where the
transition $\ell\to\ell+1$ has intensity $\ell$, and different chains
grow independently.
(In order for this to hold for $t>\tau_m$ too, we may pretend that new
items arrive also in the attic, but are ignored unless they hit an
existing chain. Similar \emph{ad hoc} modifications have to be made
after $\tau_m$ for other quantities too, in order for the
arguments below to be valid; a simple 
possibility is to redefine $\nuk(t)$ for $t>\tau_m$ so that
\eqref{b1} holds and redefine the processes $Z$ and $W$
in \eqref{erika} and \eqref{magnus}
to be constant for $t\ge\tau_m$.
The details do not matter, since we will later consider only
$t=\tau_n\le\tau_m$, so we will ignore them.)

The growth process of each chain thus is the same as the \emph{Yule process},
or \emph{binary splitting}, a branching process where each individual after a
lifetime distributed as $\Exp(1)$ splits into two.
It is well-known, see \eg{} \cite[Section III.5]{AN}, that if we start
a Yule process with a single particle at time $0$, the number of
particles at time $t$ has the geometric distribution $\Ge(\et)$ with
mean $e^t$ and
\begin{equation}
  \label{a1}
\P(k \text{ particles}) = \et\bigpar{1-\et}^{k-1},
\qquad k\ge1.
\end{equation}

Let $C(\cT)$ be the number of chains in the hash table $\cT$, and denote their
lengths by $L_1(\cT),\dots,L_{C(\cT)}(\cT)$. 
For $\ctt$ we write $C(t)$ and $L_j(t)$.

\begin{lemma}
  \label{L1}
Let $\clt$ be the number of chains of length $\ell$ in \ctt. 
Then, for each $k$ and $t\ge0$, as \mtoo,
\begin{equation}\label{l1}
  \frac1m\sum_{l\ge k} \clt
\ucpto
\int_0^{t\bmin1}(1-s)\bigpar{1-e^{-(t-s)}}^{k-1}\dd s.
\end{equation}
\end{lemma}

\begin{proof}
Informally, we observe that in a tiny time interval $[s,s+\dd s]$,
about $m\dd s$ items arrive, and $(m-N(s))\dd s\approx m(1-s)\dd s$ 
of them create new
  chains, for $s\le1$. Of these chains, by \eqref{a1}, a proportion
  $(1-e^{-(t-s)})^{k-1}$ have grown to length at least $k$ at
  time $t$, and the result follows by integration over $s$.

To be more formal, let, for $k\ge1$, $u\ge0$ and $j\ge0$, 
$\fk{u}(j)$ be the probability
that a Yule process that starts with $j$ particles
at time 0 has reached at least $k$ particles at time $u$; this is thus
equal to the probability that a 
chain of length $j$ at some instance $s$ grows to length at least
$k$ at time $s+u$. 
We have $\fk{u}(j)=1$ for $j\ge k$ and all $u$,
$\fk{0}(j)=\ind {j\ge k}$, and, by \eqref{a1}, 
\begin{equation}\label{a2}
\fk{u}(1)=\bigpar{1-\eu}^{k-1}.
\end{equation}
For fixed $k$ and $t$, and $0\le s\le t$, let
\begin{equation*}
  X(s)
\=\sum_{\ell\ge1} C_\ell(s)\fk{t-s}(\ell)
=\sum_{j=1}^{C(s)} \fk{t-s}\bigpar{L_j(s)}.
\end{equation*}
$X(s)$ is thus the expected number (given $\cT(s)$) 
of the chains present at $s$ that
have grown to length at least $k$ at $t$.
In particular, $X(t)=\cklt$.

Since new chains, all of length 1, are created with the rate 
$(m-N(s))_+$, it follows that the process
\begin{equation}\label{a3}
  Y(s)\=X(s)-\int_0^s\fk{t-u}(1)\bigpar{m-N(u)}_+\dd u,
\qquad 0\le s\le t,
\end{equation}
is a martingale. Moreover, $Y(s)$ has a jump $\Delta Y(s)=\fk{t-s}(1)$
of size $|\Delta Y(s)|\le1$
each time a new chain is created, and $Y(s)$ is smooth with a bounded
derivate between the jumps, so $s\mapsto Y(s)$ is of finite variation, 
Hence, see \eg{} Protter \cite[II.6]{Protter}, 
the quadratic variation $[Y,Y]_t=\sum_{s\le t}\Delta Y(s)^2\le N(t)$,
and, observing $Y(0)=X(0)=0$, 
\begin{equation*}
  \E Y(t)^2 = \E [Y,Y]_t \le \E N(t) = mt.
\end{equation*}
In particular, as \mtoo, $Y(t)/m\pto0$ by Chebyshev's inequality, \ie
\begin{equation*}
\frac{X(t)}{m}-\int_0^t \fk{t-s}(1)\Bigpar{1-\frac{N(s)}{m}}_+\dd s
=\frac{Y(t)}{m}
\pto 0.
\end{equation*}
Combined with \refL{LN}, this shows
\begin{equation*}
  \frac{\cklt}{m}
=
  \frac{X(t)}{m}
\pto
\int_0^t \fk{t-s}(1)\bigpar{1-s}_+\dd s,
\end{equation*}
which by \eqref{a2} proves \eqref{l1} for fixed $t\ge0$.
Convergence ucp follows by \refL{LU}.
\end{proof}

\begin{lemma}\label{LK} 
\begin{thmenumerate}
\item
For every $t\ge0$ and $r>0$, there exists a constant $K(t,r)$, not
  depending on $m$, such that
  \begin{equation*}
\E\sum_{\ell=1}^\infty \ell^r \clt 	
=
\E\sum_{j=1}^{C(t)} L_j(t)^r 
\le  K(t,r) m.
  \end{equation*}
\item
For every $r>0$, there exists a constant $K(r)$, not
  depending on $m$ or $n$, such that
  \begin{equation*}
\E\sum_{j=1}^{C(\ctmn)} L_j(\ctmn)^r 
\le  K(r) n.
  \end{equation*}
\end{thmenumerate}
\end{lemma}

\begin{proof}
\pfitem{i}
Since $Y(s)$ in \eqref{a3} is a martingale with $Y(0)=0$, 
we have $\E Y(t)=0$ and
\begin{align*}
\E C_k(t) 
&\le 
\E X(t)	 
= 
\E\intot \fk{t-u}(1)\bigpar{m-N(u)}_+\dd u
\\
&\le t\fk{t}(1)m=t(1-e^{-t})^{k-1}m.
  \end{align*}
Hence, if $a<(1-e^{-t})\qi$,
  \begin{align}\label{emma}
\E\sum_{\ell=1}^\infty a^\ell\clt
\le amt\sum_{\ell=1}^\infty\bigpar{a (\eet)}^{\ell-1}
=\frac{at}{1-a(\eet)}m
.
  \end{align}
Taking \eg{} $a=1+\et$, the result follows, 
since 
$\sup_\ell\ell^r/a^\ell<\infty$.

\pfitem{ii}
Since $\sum_j L_j(t)^r$ is increasing and $\cT(\tau_n)=\ctmn$ is
independent of $\tau_n$, we have for every $t>0$ and $a\ge1$
\begin{equation*}
\E\sum_j a^{L_j(t)}
\ge
\E\Bigpar{\sum_j a^{L_j(\tau_n)}\ind{\tau_n\le t}}
=
\E\Bigpar{\sum_j a^{L_j(\ctmn)}}\P(\tau_n\le t).
\end{equation*}
Choose $t\=2n/m\le2$ and $a\=1+e^{-2}$. Then, by \eqref{emma},
\begin{equation*}
\E\sum_j a^{L_j(t)}
\le e^4atm=2e^2(e^2+1)n.  
\end{equation*}
Moreover, $N(t)\sim\Po(mt)=\Po(2n)$, and thus
\begin{equation*}
  \P(\tau_n\le t)
=
\P\bigpar{N(t)\ge n}
=
\P\bigpar{\Po(2n)\ge n}
\to1,
\qquad
\text{as \ntoo};
\end{equation*}
hence, for some constant $c>0$ and all $n\ge1$, $\P(\tau_n\le
t)\ge c$. Consequently,
\begin{equation*}
\E\Bigpar{\sum_j a^{L_j(\ctmn)}}\le2e^2(e^2+1)c\qi n,  
\end{equation*}
and the result follows.
\end{proof}

\begin{remark}
  The collection of chain lengths evolves as a generalized P\'olya urn
  with balls of infinitely many types 0,1,\dots; we regard each empty cell as a
  ball of type 0 and each cell in a chain of length $\ell$ as a ball
  of type $\ell$. The dynamics of the urn thus is that if a ball of
  type 0 is drawn, it is removed and replaced by a ball of type 1; if a ball of
  type $\ell\ge1$ is drawn, $\ell$ balls of type $\ell$ are removed
  together with one ball of type 0, and $\ell+1$ balls of type
  $\ell+1$ are added. We start with $n$ balls of type 0.
We will, however, not use this urn representation.
\end{remark}
 
\subsection{\U}
In an unsuccessful search starting at address $j$ in a hash table
$\cT$, the number $\dujt$ of searched occupied cells is 0 if the cell
$j$ is empty; otherwise the cell belongs to a chain, and $\duj$ equals
$1$ $+$ the number of cells in the chain after $j$. 

Hence, $n_0^\U(\cT)$ is the number of empty cells in $\cT$, and 
for $\ctt$, 
\begin{equation}
  \label{nu0}
n_0^\U(t)=\bigpar{m-N(t)}_+.
\end{equation}
By \refL{LN} thus
\begin{equation}
  \label{nu00}
m\qi n_0^\U(t)\ucpto p_t^\U(0)\=(1-t)_+.
\end{equation}
(We see also that $n_0^\U(\ctmn)=m-n$, directly proving the case $k=0$
in Theorems \ref{T1}(i) and \ref{T2}(i).)

For $k\ge1$, there is exactly one cell with $\duj$ in each chain of
length $\ell\ge k$, and thus, for $\ctt$,
\begin{equation}
  \label{b1}
n_k^\U(t)=\sum_{\ell\ge k}\clt.
\end{equation}
Consequently, \refL{L1} yields, for $k\ge1$, using $u=t-s$,
\begin{equation}\label{l1a}
  \begin{split}
\frac1m n_k^\U(t)
\ucpto
 p_t^\U(k)
&\=
\int_0^{t\bmin1}(1-s)\bigpar{1-e^{-(t-s)}}^{k-1}\dd s
\\&
=
\int_{(t-1)_+}^{t}(1-t+u)\bigpar{1-e^{-u}}^{k-1}\dd u.	
  \end{split}
\end{equation}
\refT{T2}(i) follows by \refL{Ltau}.

\subsection{\li}
For the standard (late-insertion) method \li, when a new item arrives
with a hash address $j$, the insertion algorithm begins with an
unsuccessful search for the item (followed by  finding an empty cell).
The displacement of the new item is thus the same as the number $\duj$
for an unsuccessful search starting at $j$; note that for \li, the
displacement never changes after the item is inserted.
Consequently, for \ctt,
new items with displacement $k$ are created at the rate
$\nuk(t)$, and
\begin{equation}\label{erika}
  Z(t)\=\nlk(t)-\int_0^t\nuk(s)\dd s
\end{equation}
is a martingale. The jumps are all 1, and we have, as for $Y$ above,
$[Z,Z]_t\le N(t)$, and $m\qi Z(t)\pto 0$. Moreover, by Doob's 
inequality, see \eg{} \cite[p.\ 11]{Protter},
\begin{equation*}
\E\bigpar{\sup_{s\le t} Z(s)^2}
\le 4 \E Z(t)^2 = 4\E [Z,Z]_t \le 4\E N(t) = 4mt,
\end{equation*}
and hence
$m\qi Z(t)\ucpto 0$.
Consequently, by \eqref{erika} and \eqref{l1a},
\begin{equation*}
m\qi\nlk(t) \ucpto \intot p^\U_s(k)\dd s.
\end{equation*}
For $\ga>0$ we multiply by $m/n\to\gai$ and find by \refL{Ltau}
\begin{equation}\label{sofie}
n\qi\nlk(\tau_n) \pto p^\li_\ga(k)
\=\gai\intoa p^\U_s(k)\dd s
\end{equation}
as asserted in \refT{T2}(ii).
Explicitly we have, for $0<\ga\le1$,
\begin{equation*}
  \pal(0)=\gai\intoa(1-s)\dd s=1-\ga/2
\end{equation*}
and, for $k\ge1$,
\begin{align*}
\pal(k)
&
=\gai\int_{s=0}^\ga\int_{t=0}^s(1-s+t)\bigpar{1-e^{-t}}^{k-1}
\dd t \dd s
\\&
=\gai\int_{t=0}^\ga\int_{s=t}^\ga (1-s+t)\bigpar{1-e^{-t}}^{k-1}
  \dd s \dd t
\\&
=\gai\int_{t=0}^\ga \bigpar{\ga-t-(\ga-t)^2/2}\bigpar{1-e^{-t}}^{k-1}
  \dd t
.
\end{align*}

For $\ga=0$, we observe that 
$\P\bigpar{d_i^\li(\ctmn)\neq0} \le (i-1)/m$, and thus
\begin{equation*}
  \E|n-n_0^\li(\ctmn)| \le \sum_{i=1}^n i/m\le n^2/m;
\end{equation*}
hence
$  \E|1-n_0^\li(\ctmn)/n| \le n/m\to\ga=0$.
This yields \refT{T2}(ii) for $\ga=0$ with $p^\li_0(0)=1$ and
$p^\li_0(k)=0$, $k\ge1$.

\subsection{\ei}
For the early-insertion method \ei, a new item that hashes to an empty
cell gets displacement 0, which remains unchanged for ever. Hence 
$n_0^\ei(\cT)=n_0^\li(\cT)$, and
\begin{equation*}
  n\qi n_0^\ei(\cT)=n\qi n_0^\li(\cT)\pto \pal(0)=1-\ga/2
\end{equation*}
by the preceding case.

An item $x$ that hashes to an occupied cell gets an initial
displacement 1, and this displacement increases each time a new item
hashes to one of the cells in the subchain beginning with the hash
address of $x$ and ending just before $x$; the number of such cells is
the displacement, and thus the displacement grows according to the
same Yule process as the chains.
Fix $t$ and $k$, let $\fk{u}(j)$ be as in the proof of \refL{L1}, and
let, cf.\ \eqref{a3},
\begin{equation}\label{magnus}
  W(s)\=\sum_{j\ge1} n^\ei_j(s)\fk{t-s}(j)
-\int_0^s\fk{t-u}(1)N(u)\dd u.
\qquad 0\le s\le t,
\end{equation}
Again, this is a martingale. This time, however, the jumps may be
larger than 1, since more than one item can get its displacement
increased when a new item is inserted.
Clearly, the jump $\Delta W$ when a new item is inserted is at most the
length of the chain where the new item was inserted, since only the
items in this chain can have their displacements changed.
Since $[W,W]_t$ equals the sum of the squares of all jumps up to $t$,
it is at most the sum of the squares of the lengths of all chains
that have existed during the process. A chain in \ctt{} of length
$\ell$ is formed by $\ell$ insertions, and their contribution to the
latter sum is
$\sum_1^\ell k^2\le\ell^3$. Hence,
\begin{equation*}
  [W,W]_t =\sum_{s\le t}|\Delta W(s)|^2 
\le \sum_{j=1}^{C(t)} L_j(t)^3
=\sum_{\ell=1}^\infty\ell^3 \clt,
\end{equation*}
and \refL{LK} shows that
\begin{equation*}
\E W(t)^2=\E  [W,W]_t 
\le K(t,3) m.
\end{equation*}
Hence, as above,
$m\qi W(t)\pto0$, which together with \refL{LN} and \eqref{a2}
yields
\begin{equation*}
m\qi\sum_{j\ge k} n^\ei_j(t) 
 -\intot \bigpar{1-e^{-(t-u)}}^{k-1} u\dd u
\pto0.
\end{equation*}
By \refL{LU}, thus, with $s=t-u$,
\begin{equation*}
m\qi\sum_{j\ge k} n^\ei_j(t) 
\ucpto
\intot
  \bigpar{1-e^{-s}}^{k-1}(t-s)\dd s.
\end{equation*}
Replacing $k$ by $k+1$ and subtracting, we find
\begin{equation*}
m\qi n^\ei_k(t) 
\ucpto
\intot
  \bigpar{1-e^{-s}}^{k-1}e^{-s}(t-s)\dd s,
\end{equation*}
and \refT{T2}(iii) follows by \refL{Ltau} when $\ga>0$.

If $\ga=0$ we have 
$n_0^\ei(\ctmn)/n=n_0^\li(\ctmn)/n\pto1$ by the case \li,
and thus also $n_k^\ei(\ctmn)\pto0$ for $k\ge1$.

This completes the proof of \refT{T2}.

\begin{proof}[Proof of \refT{T1}]
\refT{T1}(i)--(iii) follow by taking expectations in \refT{T2} 
using dominated convergence.  

We verify directly that $\set{\pax(k)}_{k=0}^\infty$ is a probability
distribution by summing.
For \U{} we have, by \eqref{nu00} and \eqref{l1a} and summing the
geometrical series $\sum_1^\infty\bigpar{1-e^{-s}}^{k-1}=e^s$,
\begin{equation*}
  \sum_{k=0}^\infty\pau(k) 
=1-\ga+\intoa(1-\ga+s)e^s\dd s
=1-\ga+[(s-\ga)e^s]_0^\ga
=1.
\end{equation*}

For \li{} and \ei{}, the case $\ga=0$ is trivial. If $0<\ga\le1$ we
have by \eqref{sofie}
\begin{equation*}
  \sum_{k=0}^\infty\pal(k) 
=\gai\intoa\sum_{k=0}^\infty p_t^\U(k)\dd t
=\gai\intoa\dd t
=1
\end{equation*}
and, by the definition in \refT{T1},
\begin{equation*}
  \sum_{k=0}^\infty\pae(k) 
=1-\ga+\gai\intoa(\ga-t)\dd t
=1.
\end{equation*}

Finally, note that each chain of length $\ell$ contributes (for \li,
\ei{} and \U) $\ell$ displacements which are all at most
$\ell$. Hence, for $r>0$,
\begin{equation*}
  \E\bigpar{d^\Xi(\ctmn)^r\mid\ctmn}
\le \sum_j \frac{L_j}{n}L_j^r
=\frac1n \sum_j L_j^{r+1}
\end{equation*}
(also for \U), and thus, using \refL{LK}
\begin{equation*}
  \E\bigpar{d^\Xi(\ctmn)^r}
\le n\qi \E \sum_j L_j^{r+1} \le K(r+1).
\end{equation*}
Replacing $r$ by $r+1$ we see that the family $d^\Xi(\ctmn)^{r}$ is
uniformly integrable, and thus \eqref{t1} implies 
$\E (d^\Xi(\ctmn))^{r}\to\E \xpar{\dax}^r$.
\end{proof}

\begin{proof}[Proof of \refC{C1}]
It only remains to compute the expectation of $\dax$.
\begin{align*}
\E\dau
&=
\sum_{k=1}^\infty k\pau(k)
=\intoa(1-\ga+t)\sum_{k=1}^\infty k\bigpar{\eet}^{k-1}\dd t
\\
&=\intoa(1-\ga+t)\bigpar{1-(\eet)}^{-2}\dd t
=\intoa(1-\ga+t)e^{2t}\dd t
\\
&=\Bigsqpar{\Bigpar{\frac t2-\frac\ga2+\frac14}e^{2t}}_0^\ga
=\frac14 e^{2\ga}+\frac\ga2-\frac14,
\end{align*}
and, similarly, (for $\ga>0$)
\begin{align*}
\E\dal
&=
\gaix\intoa\Bigpar{\ga-t-\frac{(\ga-t)^2}2}e^{2t}\dd t
\\
&=
\gaix\Bigsqpar{\Bigpar{-\frac{(\ga-t)^2}4+\frac{\ga-t}{4}+\frac18}e^{2t}}_0^\ga
=\frac1{8\ga}e^{2\ga}+\frac\ga4-\frac14-\frac1{8\ga},
\\
\E\dae
&=
\gaix\intoa (\ga-t)e^t\dd t
=\gai\bigsqpar{\bigpar{\ga-t+1}e^{t}}_0^\ga
=\gai\bigpar{e^\ga-1-\ga}.
\end{align*}
\end{proof}

\begin{proof}[Proof of \refC{C2}]
  Similar to \refC{C1}; we omit the details.
(We used \maple{} to perform the integrations.)
For the final part, note that
\begin{align*}
\Var(\tdau)
&=\E((\dau)^2)+1-\ga-\bigpar{\E(\dau)+1-\ga}^2
\\&
=
\Var(\dau)-2(1-\ga)\E(\dau)+\ga-\ga^2.
\end{align*}
\vskip-\baselineskip
\end{proof}

\begin{proof}[Proof of \refT{T3}]
An immediate consequence of Theorems \refand{T1}{T2} and the following
general probabilistic lemma
(with $X=\dxtmn^r$, $Y=\ctmn$, $Z=\bigpar{\dax}^r$).
\end{proof}

\begin{lemma}
Let $X_n$, $Y_n$ and $Z$ be random variables (with $X_n$ and $Y_n$
defined on the same probability space, and where
$Y_n$ may take values in any measure space),
such that $X_n\ge0$ and $Z\ge0$,
and, for every real $x$, as \ntoo,
\begin{equation}
  \label{julie}
\P(X_n\le x\mid Y_n)\pto \P(Z\le x).
\end{equation}
Suppose further $\E X_n\to \E Z$.
Then $\E(X_n\mid Y_n)\pto \E Z$.
\end{lemma}

\begin{proof}
  Note first that, for every real $x$, by \eqref{julie} and dominated
  convergence, 
\begin{equation*}
\P(X_n\le x)
=\E\bigpar{\P(X_n\le x\mid Y_n)}\to \P(Z\le x),
\end{equation*}
and thus $X\dto Z$.

For any fixed $K>0$ we thus have $X_n\bmin K\dto Z\bmin K$ and thus,
by dominated convergence again,
$\E(X_n\bmin K)\to \E(Z\bmin K)$.
Hence also
\begin{equation}\label{manne}
  \begin{split}
\E(X_n-K)_+
&
=\E(X_n-X_n\bmin K)
=\E(X_n)-\E(X_n\bmin K)
\\&
\to \E(Z)-\E(Z\bmin K)
=\E(Z-K)_+.	
  \end{split}
\end{equation}

Moreover,
\begin{align*}
\bigl|\E(X_n\bmin K\mid Y_n) &- \E(Z\bmin K\mid Y_n)\bigr|
\\&
=\biggl|\int_0^K \P(X_n>x\mid Y_n)\dd x - \int_0^K \P(Z>x)\dd x\biggr|
\\&
\le \int_0^K\bigl| \P(X_n>x\mid Y_n) - \P(Z>x)\bigr|\dd x
\end{align*}
and thus, by \eqref{julie} and, yet again, dominated convergence (twice),
\begin{multline}\label{david}
\E|\E(X_n\bmin K\mid Y_n)- \E(Z\bmin K\mid Y_n)|
\\
\le \int_0^K\E\bigl| \P(X_n>x\mid Y_n) - \P(Z>x)\bigr|\dd x
\to0.
\end{multline}

By the triangle inequality and $X=X\bmin K+(X-K)_+$,
\begin{multline*}
|\E(X_n\mid Y_n)- \E(Z)|
\le
\E((X_n-K)_+\mid Y_n)
\\
+|\E(X_n\bmin K\mid Y_n)- \E(Z\bmin K)|
+ \E(Z-K)_+ .
\end{multline*}
Taking expectations, we see by \eqref{david} and \eqref{manne} that 
\begin{align*}
\limsup_\ntoo\E|\E(X_n\mid Y_n)- \E(Z)|
&\le 
\limsup_\ntoo \E(X_n-K)_+ + \E(Z-K)_+
\\&
= 2 \E(Z-K)_+
.
\end{align*}
Since $K$ is arbitrary, we see by letting $K\to\infty$ that the \lhs{}
is 0, \ie{} $\E|\E(X_n\mid Y_n)- \E(Z)|\to0$.
\end{proof}

\newcommand\AAP{\emph{Adv. Appl. Probab.} }
\newcommand\JAP{\emph{J. Appl. Probab.} }
\newcommand\JAMS{\emph{J. \AMS} }
\newcommand\MAMS{\emph{Memoirs \AMS} }
\newcommand\PAMS{\emph{Proc. \AMS} }
\newcommand\TAMS{\emph{Trans. \AMS} }
\newcommand\AnnMS{\emph{Ann. Math. Statist.} }
\newcommand\AnnPr{\emph{Ann. Probab.} }
\newcommand\CPC{\emph{Combin. Probab. Comput.} }
\newcommand\JMAA{\emph{J. Math. Anal. Appl.} }
\newcommand\RSA{\emph{Random Struct. Alg.} }
\newcommand\ZW{\emph{Z. Wahrsch. Verw. Gebiete} }
\newcommand\DMTCS{\jour{Discr. Math. Theor. Comput. Sci.} }

\newcommand\AMS{Amer. Math. Soc.}
\newcommand\Springer{Springer-Verlag}
\newcommand\Wiley{Wiley}

\newcommand\vol{\textbf}
\newcommand\jour{\emph}
\newcommand\book{\emph}
\newcommand\inbook{\emph}
\def\no#1{\unskip:#1} 

\newcommand\webcite[1]{\hfil\penalty0\texttt{\def~{\~{}}#1}\hfill\hfill}
\newcommand\webcitesvante{\webcite{http://www.math.uu.se/\~{}svante/papers/}}

\def\nobibitem#1\par{}

\end{document}